\documentclass[reqno]{amsart}
\usepackage{amssymb}
\usepackage[all]{xy}
\newtheorem{thm}{Theorem}[section]

\newtheorem{cor}[thm]{Corollary}

\def\R{{\mathbb R}}
\def\Z{{\mathbb Z}}
\def\sk{{\Sigma_k}}
\def\codim{\mathop{\rm codim}}
\def\op{{\overline\pitchfork}}
\def\wj{{\widehat\jmath}}
\def\wJ{{\widehat J}}
\def\wF{{\widehat{\op(F)}}}
\def\wW{{\widehat{W}}}
\def\Vk{{V^{(k)}}}
\def\Mk{{M^{(k)}}}
\def\DV{{\Delta^2_V}}
\def\DM{{\Delta^k_M}}
\def\fk{{f^{(k)}}}
\def\bfk{{b(f_1\times\dots\times f_k)}}
\def\ofk{{\overline\pitchfork(f_1\times\dots\times f_k)}}
\def\k{{(k)}}

\begin{document}
\title{Multiple points of immersions}
\author{Konstantin Salikhov}
\address{Department of Mathematics, University of Maryland,
College Park, MD 20742}
\email{salikhov@math.umd.edu}
\begin{abstract}
Given smooth manifolds $V^n$ and $M^m$, an integer $k>1$, and an
immersion $f:V\looparrowright M$, we have constructed
an obstruction for existence of regular homotopy of $f$ to an immersion
$f':V\looparrowright M$ without $k$-fold points. It takes values in certain
bordism group, and for $(k+1)(n+1)\leq km$ turns out to be complete.
As a byproduct, under certain dimensional restrictions we also
constructed a complete obstruction for eliminating by regular homotopy
the points of common intersection of several immersions
$f_1:V_1\looparrowright M,\dots,f_k:V_k\looparrowright M$.
\end{abstract}
\maketitle
\section{Introduction}
Let $V^n$ and $M^m$ be smooth manifolds without boundary, and $V^n$ be compact.
We will call an immersion $V\looparrowright M$ by {\it $k$-immersion}, if
it has no $k$-to-$1$ points. Then $2$-immersions $V\looparrowright M$ are
exactly embeddings $V\hookrightarrow M$. Consider the following question.
Given manifolds $V^n$, $M^m$ and integer $k>1$, do they admit a
$k$-immersion $V\looparrowright M$? Note that if $(k+1)n<km$, then
any generic immersion $V\looparrowright M$ is a $(k+1)$-immersion.
The problem of existence of an immersion $V\looparrowright M$
was solved by Hirsch \cite{Hir}. He proved that regular homotopy classes
of immersions $V\looparrowright M$ are in 1-to-1 correspondence with classes
of linear monomorphisms of tangent bundles of $\tau V\to\tau M$.
So, the natural problem would be to find out if a given regular homotopy
class of immersions $V^n\looparrowright M^m$ contains a $k$-immersion,
provided $(k+1)n<km$. Without loss of generality we can assume
that $M$ is connected.

For the case $k=2$, under a little bit stronger restrictions $3(n+1)\leq2m$,
in 1962 Haefliger gave complete answer whenever it is possible \cite{Hae}.
Haefliger's method generalizes the original purely geometric idea of Whitney
\cite{Wh} of eliminating double points of immersion.
A decade later, in 1974, Hatcher and Quinn \cite{HQ} rewrote the the
Haefliger's reasonings in the "right" language of bordism theory. In 1982
Sz\"ucs \cite{Sz 82} gave a new proof of Haefliger's theorem in the special
case $M^m=\R^m$ using the ideas of theory of singularities. Sz\"ucs informed
us that he was trying to apply his method for $k>2$. In the present paper we
generalize the ideas from \cite{HQ} and in the range $(k+1)(n+1)\leq km$
give the complete answer whenever a given regular homotopy class
of immersions $V^n\looparrowright M^m$ contains a $k$-immersion.

For a generic immersion $f:V^n\looparrowright M^m$, its set of $k$-fold
points is itself an immersed manifold $\op(\fk)\looparrowright M$.
This immersed manifold, together with some additional structure in normal
bundle of this immersion, define an element $b(\fk)$ in certain bordism
group $\Omega=\Omega_{kn-(k-1)m}(E_k;\xi(f))$. Both the group $\Omega$ and
the element $b(\fk)\in\Omega$ depend only on the regular
homotopy class of~$f$ (see Theorem~\ref{egroup}). If $(k+1)(n+1)\leq km$,
then the converse is true (Theorem~\ref{center}). Namely, if $N$ is
another $(kn-(k-1)m)$-dimensional manifold, which represents the same element
$[N]=b(\fk)\in\Omega$, then there exists an immersion $f_1$, regularly
homotopic to $f$, such that $\op(f_1^\k)=N$.

Suppose $f:V^{(k-1)r}\looparrowright M^{kr}$ is a generic immersion.
Then the manifold of $k$-fold points of $f$ is zero-dimensional. If $V,M$
are oriented and $r$ is even, to each $k$-fold point $x\in M$ of $f$ we can
attach a sign. Since $V,M$ are oriented, the normal bundle of the immersion
$f$ is oriented. Since $r$ is even, the orientation of the normal spaces to
$k$ "sheets" of $V$, intersecting at $x$, gives an orientation of the tangent
space $\tau_x M$. We put sign "$+$" if this orientation coincides with the
orientation of $M$, and "$-$" otherwise. Define $I(f)\in\Z$ as the number
of $k$-fold points, counted with signs. If either $V^{(k-1)r}$ or $M^{kr}$
is not orientable, or $r$ is odd, define $I(f)\in\Z_2$ as the number of
$k$-fold points modulo 2.
\begin{cor}\label{general0}
Let $V^{(k-1)r},M^{kr}$ be connected smooth manifolds without boundary,
$V$ is compact, $M$ is simply connected, and $r>k$. Then a generic immersion
$f:V\looparrowright M$ is regularly homotopic to a $k$-immersion iff $I(f)=0$.
\end{cor}
\begin{cor}\label{odd0}
Let $V^{(k-1)r}$ be a smooth compact orientable manifold without boundary,
and $r>k$ are odd integers. Then any immersion $f:V\looparrowright\R^{kr}$
is regularly homotopic to a $k$-immersion.
\end{cor}
A map $f:V\to M$ is called a {\it "topological immersion"}, if for any
point $x\in V$ there exists an open neighborhood $U_x\ni x$ such that
for any $y_1\neq y_2\in U_x$ we have $f(y_1)\neq f(y_2)$.
\begin{cor}\label{topimm}
Let $V^n,M^m$ be smooth manifolds without boundary, $V^n$ be
compact, $f:V\looparrowright M$ be an immersion, and $(k+1)(n+1)\leq km$.
Then $f$ is regularly homotopic to a $k$-immersion iff there exists a
"regular homotopy" of $f$ through topological immersions to a
map without $k$-fold points.
\end{cor}
The main problem of constructing the desired regular homotopy virtually
splits into two parts: geometrical and combinatorial. The first part deals
with construction of this regular homotopy locally. The second part, which
in general turns out to be harder, deals with the fact that the natural
$k$-fold covering map $\widetilde{\op(\fk)}\to\op(\fk)$ can be non-trivial.
In the section~\ref{simp} we restrict ourselves to the case, then the second
part is trivial. Namely, we try to figure out whenever immersions of $k$
different manifolds $f_1:V_1\looparrowright M,\dots,f_k:V_k\looparrowright M$
can be regularly homotoped to immersions without common intersection.
\begin{cor}\label{multhom}
Let $V_i^{n_i}$, $i=1\dots k$ and $M^m$ be smooth $C^\infty$ manifolds
without boundary, all $V_i$ are compact, $2(p+1)<n_i<m-p-k$ for all
$i=1\dots k$, and $2(p+k)<m$, where $p=n_1+\dots+n_k-(k-1)m$. Suppose
$f_i:V_i\looparrowright M$ are smooth immersions. Then the immersions
$f_1,\dots,f_k$ are regularly homotopic to immersions without common
intersection iff $f_1,\dots,f_k$ can be continuously homotoped to
continuous maps with empty common intersection.
\end{cor}
\begin{cor}\label{multsh}
Let $V_i^{n_i}$, $i=1\dots k$ and $M^m$ be smooth $C^\infty$ manifolds
without boundary, all $V_i$ are compact, $2(p+1)<n_i<m-p-k$ for all
$i=1\dots k$, and $2(p+k)<m$, where $p=n_1+\dots+n_k-(k-1)m$. Suppose
$p\in\{4,5,12\}$, all $V_i$ are $p$-connected, and $M$ is $(p+1)$-connected.
Then any immersions $f_i:V_i\looparrowright M$ can be regularly homotoped
to immersions without common intersection.
\end{cor}
\section{Simple case: mutual intersections of $k$ manifolds}\label{simp}
Let $T$ be $(k-1)$-dimensional simplex, $\sigma_1,\dots,\sigma_k$ be its
vertices and $\sigma$ be its barycenter.
For any topological spaces $V_1,V_2,\dots,V_k,M$ and continuous maps
$f_1:V_1\to M,\dots,f_k:V_k\to M$, consider the space
\begin{equation*}
E=E(f_1,\dots,f_k)=\{(x_1,\dots,x_k,\theta)\mid x_i\in V_i\text{ and }
\theta:T\to M\text{ such that }\theta(\sigma_i)=f_i(x_i)\}
\end{equation*}
Then there are obvious projections $\pi_i:E\to V_i$ and
$\pi:E\to\{\theta:T\to M\}$. The map $\Phi:E\times T\to M$, defined by
$(e,t)\mapsto\pi(e)(t)$, where $e\in E$ and $t\in T$,
makes the diagram below homotopy commutative
\begin{equation*}
\xymatrix{
&&E\ar[dll]_{\pi_1}\ar[dl]^{\pi_2}\ar[dr]^{\pi_k}\\
V_1\ar[drr]_{f_1}&V_2\ar[dr]^{f_2}&\dots&V_k\ar[dl]^{f_k}\\
&&M
}
\end{equation*}
The space $E$ is universal in the following sense. If we are given a
space $X$, maps $h_i:X\to V_i$, and a homotopy $H:X\times T\to M$
such that $H(\cdot,\sigma_i)=f_ih_i$ (which connects all $f_ih_i$
with some map $H(\cdot,\sigma)$), then there is a unique map
$j:X\to E$  such that $H(x,t)=\Phi(j(x),t)$,
which is defined by $j(x)=(h_1(x),\dots,h_k(x),H(x,\cdot))$.
\subsection{Bordism group}
From now, suppose all $V_i$ and $M$ are smooth $C^\infty$ manifolds
without boundary and $V_i$ are compact. Denote by $\nu_{V_i}$
normal bundles over $V_i$, and by $\tau_M$ the tangent bundle over $M$.
Over $E$, consider the bundle
$\xi=\pi_1^*(\nu_{V_1}\oplus f_1^*\tau_M)\oplus\dots\oplus
\pi_k^*(\nu_{V_k}\oplus f_k^*\tau_M)\oplus\Phi(\cdot,\sigma)^*\nu_M$.
Define the bordism group $\Omega_p(E;\xi)$, whose objects are tuples
$(N^p,\nu_N,\gamma,\omega_N)$, where $N$ is a compact manifold without
boundary, $\nu_N$ is a normal bundle over $N$, $\gamma:N\to E$ is a map,
and $\omega_N:\nu_N\overset{\sim}\to\gamma^*\xi$ is a stable isomorphism.
Note that this group depends only on the homotopy classes of $f_1,\dots,f_k$.
By Pontryagin-Thom construction, $\Omega_p(E;\xi)=
\lim_{s\to\infty}\pi_{p+\dim\xi+s}T(\xi\oplus\varepsilon^s)$, where
$T$ denotes the Thom space, and $\varepsilon$ is the trivial bundle.

Denote by $\DM$ the diagonal $\{(x,\dots,x)\in\Mk\mid x\in M\}$.
If the map $f_1\times\dots\times f_k:V_1\times\dots\times V_k\to\Mk$ is
smooth and transversal to the diagonal $\DM$, then we
say {\it "the maps $f_1,\dots,f_k$ are transversal"} and denote the manifold
$(f_1\times\dots\times f_k)^{-1}(\DM)\subset V_1\times\dots\times V_k$
by $\ofk$. Since all $V_i$ are compact, $\ofk$ is also compact.
Note that the notion for maps $f_1,\dots,f_k$ to be transversal is generic,
and any kit of maps $f_1,\dots,f_k$ can be approximated by a transversal one.
\begin{thm}\label{group}
Suppose $V_1,\dots,V_k,M$ are smooth $C^\infty$ manifolds without boundary,
$V_1,\dots,V_k$ are compact, and the maps $f_i:V_i\to M$, $i=1\dots k$
are continuous. Then this defines a canonical element
$\bfk\in\Omega_*(E(f_1,\dots,f_k);\xi)$.
This element depends only on the homotopy classes of $f_1,\dots,f_k$.
\end{thm}
\begin{thm}\label{easy}
Let $V_i^{n_i}$, $i=1\dots k$ and $M^m$ be smooth $C^\infty$ manifolds
without boundary, all $V_i$ are compact, $2(p+1)<n_i<m-p-k$ for all
$i=1\dots k$, and $2(p+k)<m$, where $p=n_1+\dots+n_k-(k-1)m$. Suppose
$f_i:V_i\looparrowright M$ are smooth immersions. Let $\gamma:N^p\to E$ be
a singular manifold, and $\omega_N:\nu_N\overset{\sim}\to\gamma^*\xi$ be a
stable isomorphism such that $[N]=\bfk\in\Omega_p(E;\xi)$. Then there exists
a regular homotopy of the immersions $f_1,\dots,f_k$ to transversal
immersions $f_1',\dots,f_k'$ such that $\op(f_1'\times\dots\times f_k')=N$.
\end{thm}
\section{General case: $k$-fold self-intersections}
Let $V$ and $M$ be topological spaces and $f:V\to M$ be a continuous map.
The group $\sk$ of permutations on $k$ elements $\{1,\dots,k\}$ acts on
the $(k-1)$-simplex $T$ linearly, permuting its vertices
$\sigma_1,\dots,\sigma_k$. Note that the only fixed point for this action
is the barycenter $\sigma$ of $T$. Take $k$ copies of the manifold $V$ and
enumerate them by $1\dots k$ in arbitrary way $V_1,\dots,V_k$.
Consider the space
\begin{equation*}
\widehat{E_k}=\widehat{E_k}(f)=\{(x_1,\dots,x_k,\theta)\mid x_i\in V_i,
x_i\neq x_j \text{ for }i\neq j\text{, and }
\theta:T\to M\text{ such that }\theta(\sigma_i)=f(x_i)\}
\end{equation*}
Define the action of element $g\in\sk$ on $\widehat{E_k}$ by the formula
$g(x_1,\dots,x_k,\theta)=
(x_{g^{-1}(1)},\dots,x_{g^{-1}(k)},\theta\circ g^{-1})$.
Note that this action is free. Put $E_k=\widehat{E_k}/\sk$. Define the
projections $\pi_i:\widehat{E_k}\to V$ and
$\pi:\widehat{E_k}\to\{\theta:T\to M\}$ by the formulas
$\pi_i(x_1,\dots,x_k,\theta)=x_i$ and $\pi(x_1,\dots,x_k,\theta)=\theta$.
It is easy to see, the map $\Phi:\widehat{E_k}\times T\to M$, defined by
$(e,t)\mapsto\pi(e)(t)$, is $\sk$-invariant with respect to the diagonal
$\sk$-action on $\widehat{E_k}\times T$. Therefore the map $\phi:E_k\to M$ is
well-defined by the formula $[e]\mapsto\Phi(e,\sigma)$. Note that the homotopy
types of $\widehat{E_k}$ and $E_k$ depend only on the homotopy class of $f$.

The space $\widehat{E_k}$ is universal in the following sense. For any space
$X$ with a free $\sk$-action denote by $X_{(i)}$ the quotion space
$X/\Sigma^i_{k-1}$, where $\Sigma^i_{k-1}\subset\sk$ is the stabilizer of the
element $i$. Denote by $\pi_i':X\to X_{(i)}$ and by $\pi':X\to X/\sk$ the
natural projections. Denote by $\Sigma^i_{k-1}(x)$ the orbit of $x\in X$
under $\Sigma^i_{k-1}$-action. Note that
$g(\Sigma^i_{k-1}(x))=\Sigma^{g(i)}_{k-1}(g(x))$. So, for any $g\in\sk$ the map
$g:X_{(i)}\to X_{(g(i))}$ is well-defined, and it is identity if $g=1\in\sk$.
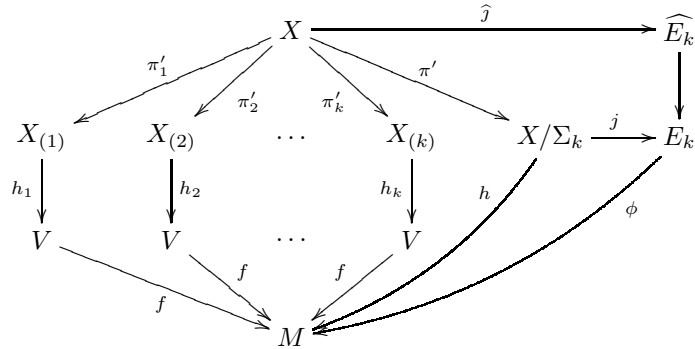
\begin{equation*}
\xymatrix{
&&X\ar[dll]_{\pi_1'}\ar[dl]^{\pi_2'}\ar[dr]_{\pi_k'}\ar[drr]^{\pi'}\ar[rrr]^{\wj}&&&\widehat{E_k}\ar[d]\\
X_{(1)}\ar[d]_{h_1}&X_{(2)}\ar[d]^{h_2}&\dots&X_\k\ar[d]_{h_k}&X/\sk\ar@/^4mm/[ddll]_<<<<<<<{h}\ar[r]^{j}&E_k\ar@/^5mm/[ddlll]^<<<<<<{\phi}\\
V\ar[drr]_{f}&V\ar[dr]^{f}&\dots&V\ar[dl]_{f}\\
&&M
}
\end{equation*}
Suppose, we are given the following data:
\begin{itemize}
\item a space $X$ with a free $\sk$-action
\item a map $h:X/\sk\to M$
\item maps $h_i:X_{(i)}\to V_i(\cong V)$ such that $h_{g(i)}\circ g=h_i$ for
any $g\in\sk$, and $h_i\pi_i'(x)\neq h_j\pi_j'(x)$ for $x\in X$ and $i\neq j$
\item a $\sk$-invariant homotopy $H:X\times T\to M$
such that $H(\cdot,\sigma)=h\pi'$ and $H(\cdot,\sigma_i)=fh_i\pi_i'$
\end{itemize}
Then there is a unique $\sk$-equivariant map $\wj:X\to\widehat{E_k}$
such that $H(x,t)=\Phi(\wj(x),t)$, which is defined by
$\wj(x)=(h_1\pi_1'(x),\dots,h_k\pi_k'(x),H(x,\cdot))$. Now we can get rid of
the ambiguity of the choice of enumeration $V_1,\dots,V_k$. Clearly, if we
will take another enumeration, which differs from the original one by
$g\in\sk$, then the "classifying map" $\wj$ will become $\wj\circ g=g\circ\wj$.
Thus, all the information about $X$ is actually equivalent to a $\sk$-orbit of
$\sk$-equivariant maps $\wj:X\to\widehat{E_k}$. Since $\sk$-actions on $X$ and
$\widehat{E_k}$ are free, this is equivalent to the map
$j:=\wj/_\sk:X/\sk\to\widehat{E_k}/\sk=E_k$.
\subsection{Bordism group}
From now, suppose $V^n$ and $M^m$ are smooth $C^\infty$ manifolds without
boundary, $V$ is compact. Let $f:V\looparrowright M$ be an immersion. Denote
by $\nu_f=\nu_f(V,M)$ the normal bundle over $V$, induced by the immersion
$f$, and by $\nu_M$ a normal bundle over $M$. Over $\widehat{E_k}$, consider
the bundle $\pi_1^*\nu_f\oplus\dots\oplus\pi_k^*\nu_f$. Define the $\sk$-action
in the total space of this bundle by $g(e,\vec{v}_1,\dots,\vec{v}_k)=
(g(e),\vec{v}_{g^{-1}(1)},\dots,\vec{v}_{g^{-1}(k)})$. Since this action
covers the $\sk$-action on $\widehat{E_k}$, this gives a well-defined vector
bundle $(\pi_1^*\nu_f\oplus\dots\oplus\pi_k^*\nu_f)/\sk$ over
$E_k=\widehat{E_k}/\sk$. Note that the structural group of the last bundle
is $O(m-n)\wr\sk$ (see \cite{AE,KS} for details). As a set,
$O(m-n)\wr\sk=O(m-n)^\k\times\sk$, and the multiplication is defined by
$(A_1,\dots,A_k,\sigma)*(B_1,\dots,B_k,\tau)=
(A_1B_{\sigma^{-1}(1)},\dots,A_kB_{\sigma^{-1}(k)},\sigma\tau)$.
Over $E_k$, consider the bundle
$\xi=\xi(f)=(\pi_1^*\nu_f\oplus\dots\oplus\pi_k^*\nu_f)/\sk\oplus\phi^*\nu_M$.
Its structural group is $O(m-n)\wr\sk\oplus O$. Denote $p=kn-(k-1)m$.
Define the bordism group $\Omega_p(E_k;\xi)$, whose objects are tuples
$(N^p,\nu_N,\gamma,\omega_N)$, where $N$ is a compact manifold without
boundary, $\nu_N$ is a normal bundle over $N$, $\gamma:N\to E_k$ is a map,
and $\omega_N:\nu_N\overset{\sim}\to\gamma^*\xi$ is a stable isomorphism.
Note that the bundle $\nu_f$, and, therefore, the group $\Omega_p(E_k;\xi(f))$,
depend only on the regular homotopy class of $f$. By Pontryagin-Thom
construction, $\Omega_p(E_k;\xi)=
\lim_{s\to\infty}\pi_{p+\dim\xi+s}T(\xi\oplus\varepsilon^s)$, where
$T$ denotes the Thom space, and $\varepsilon$ is the trivial bundle.

Denote by $\DV$ the diagonal $\{(x_1,\dots,x_k)\in\Vk\mid
x_i\neq x_j\text{ for }i\neq j\}$. Note that in the important special case
$M^m=\R^m$ the space $E_k$ reduces to $kn$-dimensional manifold
$(\Vk-\DV)/\sk$, and $\xi$ reduces to the bundle $\nu_\fk/\sk$.

We say that a $\sk$-equivariant map $F:\Vk\to\Mk$ is {\it "$k$-disjoint"},
if $F^{-1}(\DM)-\DV$ is disjoint from $\DV$. Note that for any topological
immersion $f:V\to M$, the map $\fk$ is $k$-disjoint. We say that a
$k$-disjoint map $F$ is {\it "$k$-transversal"}, if either
$F^{-1}(\DM)-\DV=\emptyset$, or the map $F|_{\Vk-\DV}$ is smooth and
transversal to~$\DM$. By equivariant transversality theorem \cite{GG},
for a generic smooth immersion $f:V\looparrowright M$, the map $\fk$ is
$k$-transversal. Denote the manifold $F^{-1}(\DM)-\DV$ by $\wF$. Since $\Vk$
is compact and $\wF$ is disjoint from $\DV$, the manifold $\wF$ is also
compact. Since $\wF$ is the intersection of $\sk$-invariant set $\Vk-\DV$
with the preimage of $\sk$-equivariant map $F:\Vk\to\Mk$ of the
$\sk$-invariant set $\DM$, it is invariant under $\sk$-action on $\Vk$.
Since the $\sk$-action on $\Vk-\DV$ is free, its restriction on $\wF$ is
also free. Denote $\op(F)=\wF/\sk$. Note that if $f:V\looparrowright M$ is
a generic immersion, then $\op(\fk)\looparrowright M$ is the locus of
its $k$-fold points.
\begin{thm}\label{egroup}
Suppose $V^n,M^m$ are smooth $C^\infty$ manifolds without boundary, $V$ is
compact, and $F:\Vk\to\Mk$ is a $k$-transversal map, $\sk$-equivariantly
homotopic to to $\fk$ for some immersion $f:V\looparrowright M$. Then
$\op(F)$ defines a canonical element $b(F)\in\Omega_p(E_k;\xi(f))$.
This element depends only on the $k$-disjoint homotopy class of $F$.
\end{thm}
Note that for an immersion $f:V\looparrowright M$, the element
$b(\fk)\in\Omega_p(E_k;\xi(f))$ is well-defined even if $\fk$ is not
$k$-transversal. Indeed, let $f_1$ be an immersion, regularly homotopic to
$f$, such that $f_1^\k$ is $k$-transversal. From the second part of
Theorem~\ref{egroup} it follows that $b(f_1^\k)$ doesn't depend on the
choice of $f_1$. Put $b(\fk)=b(f_1^\k)\in\Omega_p(E_k;\xi(f))$.
\begin{thm}\label{center}
Let $V^n,M^m$ be smooth $C^\infty$ manifolds without boundary, $V$ be
compact, $(k+1)(n+1)\leq km$, and $f:V\looparrowright M$ be an immersion.
Let $\gamma:N^p\to E_k$ be a singular manifold, and
$\omega_N:\nu_N\overset{\sim}\to\gamma^*\xi$ be a stable isomorphism
such that $[N]=b(\fk)\in\Omega_p(E_k;\xi(f))$. Then there exists a regular
homotopy of the immersion $f$ to an immersion $f_1$ such that $f_1^\k$ is
$k$-transversal and $N=\op(f_1^\k)$.
\end{thm}
\section{Proofs of Corollaries}
\begin{proof}[Proof of Corollary~\ref{general0}]
Note that $\pi_1\times\dots\times\pi_k:\widehat{E_k}\to\Vk-\DV$ is
a Serre fibration with the fiber, homotopy equivalent to
$F=\{\theta:T\to M\mid\theta(\sigma_i)=*\}$, where $*\in M$ is a base point.
In $T$, consider $1$-dimensional subcomplex $T_1$, consisting of vertices
$\sigma_1,\dots,\sigma_k$ and $(k-1)$ edges
$\sigma_1\sigma_2,\dots,\sigma_1\sigma_k$. Since $T_1$ is a deformation
retract of $T$, the space $F$ is homotopy equivalent to
$F_1=\{\theta:T_1\to M\mid\theta(\sigma_i)=*\}=(\Omega M)^{(k-1)}$.
Since $M$ is $1$-connected, then the loopspace $\Omega M$ is connected.
Since $V$ is connected and $\codim_\Vk(\DV)=n>1$, then $\Vk-\DV$ is
connected. Hence $\widehat{E_k}$ and $E_k=\widehat{E_k}/\sk$ are connected.

Since $M$ is $1$-connected, then $\phi^*\nu_M$ is always orientable. Therefore
$\xi=(\pi_1^*\nu_f\oplus\dots\oplus\pi_k^*\nu_f)/\sk\oplus\phi^*\nu_M$ is
orientable iff $(\pi_1^*\nu_f\oplus\dots\oplus\pi_k^*\nu_f)/\sk$ is orientable.
Suppose $V$ is orientable and $(m-n)$ is even. Then $\nu_f$ is orientable,
fix its orientation. Then $\pi_i^*\nu_f$ are oriented. Since $(m-n)$ is even,
the $\sk$-action on $\pi_1^*\nu_f\oplus\dots\oplus\pi_k^*\nu_f$ preserves
orientation. Hence $(\pi_1^*\nu_f\oplus\dots\oplus\pi_k^*\nu_f)/\sk$ is
orientable.

Suppose $V$ is orientable and $(m-n)$ is odd. Then $\nu_f$ is orientable,
fix its orientation. Then $\pi_i^*\nu_f$ are oriented. Choose points
$b_0,b_1\in\widehat{E_k}$ such that $g(b_0)=b_1$, where $g\in\sk$
is an odd permutation. Since $\widehat{E_k}$ is connected, there exists a
path $b:[0,1]\to\widehat{E_k}$ such that $b(0)=b_0,b(1)=b_1$. Then the
projection of $b_t$ into $E_k$ defines a loop. Since $r$ is odd
and $g$ is odd, then the restriction of
$(\pi_1^*\nu_f\oplus\dots\oplus\pi_k^*\nu_f)/\sk$ on this loop
is non-orientable.

Suppose $V$ is non-orientable. Then $\nu_f$ is non-orientable. Choose a
loop $c:[0,1]\to V$ such that the restriction of $\nu_f$ on $c_t$ in
non-orientable. Choose points $c_2,\dots,c_k\in V$ such that $c_i\neq c_j$
if $i\neq j$, and $c_i\neq c(t)$ for all $t\in[0,1]$. This defines a loop
$C:[0,1]\to\Vk-\DV$ via $C(t)=(c(t),c_2,\dots,c_k)$. Since
$\pi_1\times\dots\times\pi_k:\widehat{E_k}\to\Vk-\DV$ is a Serre fibration
with connected fiber, then there exists a loop $B:[0,1]\to\widehat{E_k}$
such that $\pi_1\times\dots\times\pi_k\circ B_t\simeq C_t$. Then the
projection of the loop $B_t$ into $E_k$ defines a loop, along which
$(\pi_1^*\nu_f\oplus\dots\oplus\pi_k^*\nu_f)/\sk$ is non-orientable.

Choose a base point $*\in E_k$. Since $E_k$
is connected, any map of a $0$-dimensional manifold $\gamma:N^0\to E_k$ can
be homotoped to the constant map $N^0\to *$. Since the orthogonal group
$O$ has two connected components, then a two-point singular manifold
$\{a_1,a_2\}\to *$ bounds iff the orientations of normal bundles
of $a_1,a_2$ at~$*$ are opposite. Suppose the bundle $\xi$ is orientable.
Fix its orientation. Then the orientation of the normal bundle of $a_1$
at~$*$ does not depend on the choice of homotopy
$\gamma\rightsquigarrow(N^0\to *)$.
Therefore the bordism class of a singular manifold $\gamma:N^0\to E_k$
is completely characterized by the number $l_1-l_2\in\Z$, where $l_1$ is
the number of points in $N^0$, whose orientation coincides with the
orientation of $\xi$, and $l_2$ --- with opposite orientation. Hence
$\Omega_0(E_k;\xi)=\Z$. If $\xi$ is non-orientable, then any singular
manifold $a\to *$ is bordant to itself with opposite orientation via a loop,
that changes the orientation of $\xi$. Therefore $\Omega_0(E_k;\xi)=\Z_2$.
It is easy to see that the
invariant $b(\fk)\in\Omega_0(E_k;\xi)$ coincides with $I(f)$. Since $r>k$,
then $(k+1)((k-1)r+1)\leq k^2r$. Then Corollary~\ref{general0} follows from
Theorems~\ref{egroup} and~\ref{center}.
\end{proof}
\begin{proof}[Proof of Corollary~\ref{odd0}]
WLOG we may assume that $V$ is connected. By Corollary~\ref{general0},
it suffices to show that the number of $k$-fold points of any generic
immersion $V\looparrowright\R^{kr}$ is even. This is exactly the result
of \cite{Sz 90}.
\end{proof}
\begin{proof}[Proof of Corollary~\ref{topimm}]
Let $f_t:V\to M,t\in[0,1]$ be a topological regular homotopy, $f_0=f$ and
$f_1$ has no $k$-to-$1$~points. Then $f_t^\k$ is a $k$-disjoint homotopy
from $\fk$ to a map $f_1^\k$ such that $(f_1^\k)^{-1}(\DM)-\DV=\emptyset$.
Then by Theorem~\ref{egroup}, the class $b(\fk)=0\in\Omega_p(E_k;\xi)$,
and Theorem~\ref{center} gives the desired regular homotopy.
\end{proof}
\begin{proof}[Proof of Corollary~\ref{multhom}] follows immediately from
Theorem~\ref{group} and Theorem~\ref{easy}.
\end{proof}
\begin{proof}[Proof of Corollary~\ref{multsh}]
As we saw above, $E\to V_1\times\dots\times V_k$ is a Serre fibration with
fiber, homotopy equivalent to $(\Omega M)^{(k-1)}$. Since all $V_i$ are
$p$-connected and $M$ is $(p+1)$-connected, then $E$ is $p$-connected. Then
any map $\gamma:N^p\to E$, representing a bordism class $\Omega_p(E;\xi)$, is
null-homotopic. Therefore $\nu_N$ is trivial, and, by Pontryagin-Thom
construction, $(N,\omega_N)$ represents an element of $p$-th stable homotopy
group of spheres. Since this group vanishes for $p\in\{4,5,12\}$, then
$(N,\omega_N)=0$. Therefore $\Omega_p(E;\xi)=0$, and by Theorem~\ref{easy},
any immersions $f_i:V_i\looparrowright M$ are regularly homotopic to
immersions without common intersection.
\end{proof}
\section{Proofs of Theorems}
\begin{proof}[Proof of Theorem~\ref{group}]
First, suppose the maps $f_1,\dots,f_k$ are transversal. Then there are
obvious projections $h_i:\ofk\to V_i$ such that all $f_ih_i$ are connected
by constant homotopy. By universal property of $E(f_1,\dots,f_k)$, there is
a canonical map $j:\ofk\to E$. Denote by $\nu(\DM,\Mk)$ the normal bundle
of $\DM$ in $\Mk$. By construction of $\ofk$, the normal bundle
$\nu(\ofk,V_1\times\dots\times V_k)$ is
$(f_1\times\dots\times f_k)^*\nu(\DM,\Mk)\cong
j^*(f_1\pi_1\times\dots\times f_k\pi_k)^*\nu(\DM,\Mk)\simeq
j^*(f_1\pi_1\times\dots\times f_k\pi_k)^*(\tau_\Mk\oplus\nu_\DM)$.
But $(f_1\pi_1\times\dots\times f_k\pi_k)\circ j=
D\circ\Phi(\cdot,\sigma)\circ j$, where $D:M\to\DM$ is the diffeomorphism
$x\mapsto(x,\dots,x)$. Then $\nu_\ofk=
\nu_{V_1\times\dots\times V_k}|_\ofk\oplus
\nu(\ofk,V_1\times\dots\times V_k)\simeq
j^*(\pi_1\times\dots\times\pi_k)^*(\nu_{V_1\times\dots\times V_k})
\oplus j^*(f_1\pi_1\times\dots\times f_k\pi_k)^*\tau_\Mk\oplus
j^*\Phi(\cdot,\sigma)^*\nu_M=j^*\xi$. Thus the map $j:\ofk\to E$ defines
an element $\bfk\in\Omega_p(E;\xi)$.

Now, suppose we have a kit of homotopies $f_{1,t},\dots,f_{k,t}$, $t\in[0,1]$,
such that the kits $f_{1,0},\dots,f_{k,0}$ and $f_{1,1},\dots,f_{k,1}$ are
transversal. Then it can be approximated by a new one such that
$F=f_{1,\cdot}\times\dots\times f_{k,\cdot}:
V_1\times\dots\times V_k\times[0,1]\to\Mk$
is transversal to $\DM$ and is left fixed on both ends $t=0,1$ (where we
already have transversality). Let the manifold $W=F^{-1}(\DM)$. Then
$\partial W=\op(f_{1,0}\times\dots\times f_{k,0})\cup
\op(f_{1,1}\times\dots\times f_{k,1})$.

Denote by $t:V_1\times\dots\times V_k\times[0,1]\to[0,1]$ the projection
on the last factor. Here we again have obvious projections
$\tilde h_i:W\to V_i$ such that $f_{i,0}\tilde h_i$ are homotopic to
$D^{-1}\circ F=f_{i,t(\cdot)}(h_i(\cdot))$ by homotopy
$x\mapsto f_{i,t(x)t'}(h_i(x))$, $t'\in[0,1]$. This gives a homotopy
$H:W\times T^1\to M$, where $T^1$ is the $1$-dimensional subcomplex of $T$,
consisting of vertices $\sigma,\sigma_1,\dots,\sigma_k$ and straight
intervals $\sigma\sigma_1,\dots,\sigma\sigma_k$. Since there exists a
canonical conical retraction of $T$ onto $T^1$, there is a canonical
extension of $H$ to a homotopy $H:W\times T\to M$. By universal property
of $E$, this gives a canonical map $J:W\to E$.

The normal bundle $\nu(W,V_1\times\dots\times V_k\times[0,1])$ is
$F^*\nu(\DM,\Mk)\cong(F\circ(\pi_1\times\dots\times\pi_k\times t)\circ J)^*
\nu(\DM,\Mk)\simeq J^*(\pi_1\times\dots\times\pi_k\times t)^*F^*
(\tau_\Mk\oplus\nu_\DM)$. But $F\circ(\pi_1\times\dots\times\pi_k\times t)
\circ J= D\circ\Phi(\cdot,\sigma)\circ J$.
Note that $(\pi_1\times\dots\times\pi_k\times t\circ J)^*
(\nu_{V_1\times\dots\times V_k\times[0,1]})\simeq
(\pi_1\times\dots\times\pi_k\circ J)^*(\nu_{V_1\times\dots\times V_k})$
and $(\pi_1\times\dots\times\pi_k\times t\circ J)^*F^*
\tau_\Mk\simeq(\pi_1\times\dots\times\pi_k\circ J)^*
(f_{1,0}\times\dots\times f_{k,0})^*\tau_\Mk$. Then
$\nu_W=\nu_{V_1\times\dots\times V_k\times[0,1]}|_W\oplus
\nu(W,V_1\times\dots\times V_k\times[0,1])\simeq
J^*(\pi_1\times\dots\times\pi_k\times t)^*
(\nu_{V_1\times\dots\times V_k\times[0,1]})\oplus
J^*(\pi_1\times\dots\times\pi_k\times t)^*F^*\tau_\Mk
\oplus J^*\Phi(\cdot,\sigma)^*\nu_M\simeq J^*\xi$.
So, $J:W\to E$ gives a bordism between $b(f_{1,0}\times\dots\times f_{k,0})$
and $b(f_{1,1}\times\dots\times f_{k,1})$ in $\Omega_p(E(f_1,\dots,f_k);\xi)$.

Finally, recall that any kit of maps $f_1,\dots,f_k$ can be approximated
by a transversal kit $f_1',\dots,f_k'$, homotopic to $f_1,\dots,f_k$. Define
$\bfk:=b(f_1'\times\dots\times f_k')$. The second half of the theorem shows,
that it does not depend on the choice of the transversal approximation
$f_1',\dots,f_k'$.
\end{proof}
\begin{proof}[Proof of Theorem~\ref{easy}]
First, make a small regular perturbation of the immersions $f_1,\dots,f_k$
to put them into general position. By Theorem~\ref{group} it will not
change the class $\bfk\in\Omega_p(E;\xi)$. Now the immersions $f_1,\dots,f_k$
are transversal, and the manifold $\ofk$ is defined. Since $2p<n_i<m$ for
all $i=1\dots k$, by general position we may assume that $\ofk$ is embedded
into each $V_i$ and into $M$. Since $n_i<m-p$, we have
$2n_i-m+(n_1+\dots+\widehat{n_i}+\dots+n_k-(k-2)m)<m$
for all $i=1\dots k$. This means that by general position we may assume
that $(2n_i-m)$-dimensional manifold of self-intersections of
$f_i:V_i\looparrowright M$ is disjoint from
$(n_1+\dots+\widehat{n_i}+\dots+n_k-(k-2)m)$-dimensional immersed manifold
of mutual intersections of $f(V_1),\dots,\widehat{f_i(V_i)},\dots,f_k(V_k)$
in $M^m$.

Let $J:W^{p+1}\to E$ be the bordism between $\ofk$ and $N$ in
$\Omega_p(E;\xi)$. Define the map $H:W\times T\to M$ by the formula
$H(w,t)=\Phi(J(w),t)$. Note that $H(W,\sigma_i)\subset f_i(V_i)$
and $H|_\ofk$ is a constant homotopy.

Since for all $i=1\dots k$ we have $2(p+1)<n_i<m-(p+1)$, then
$2\dim W<\dim(V_i)$, the dimension $\dim W+\dim(\text{self-intersections
of $V_i$ in $M$})<\dim(V_i)$, and $\dim W+\dim(\text{intersections of $V_i$
with $V_j$})<\dim(V_i)$. Again, applying general position argument, we may
$C^0$-perturb the map $H:W\times T\to M$ (leaving it fixed on $\ofk\times T$)
such that $H|_{W\times\sigma_i}$ will be smooth embeddings
$W\times\sigma_i\hookrightarrow f_i(V_i)$, disjoint from self-intersections
of $f_i:V_i\looparrowright M$, and intersections with over $V_j$, distinct
from $\ofk$.

Finally, since $2(p+k)<m$ and $n_i<m-p-k$ for all $i=1\dots k$, we have
$2\dim(W\times T)<\dim(M)$ and $\dim(W\times T)+\dim(V_i)<\dim(M)$. This means
that by general position we can $C^0$-perturb the map $H$ (leaving it fixed
on $\ofk\times T$ and $W\times\sigma_i$) such that
$H:(W-\ofk)\times T\hookrightarrow M$ will be a smooth embedding with
$H(W\times T-\ofk\times T-W\times\{\sigma_1,\dots,\sigma_k\}))\cap
(\bigcup_{i=1\dots k}f_i(V_i))=\emptyset$.

Denote by $W^+$ and $T^+$ open manifolds without boundary, which are obtained
from $W$ and $T$ by attaching a small collar neighborhood of the boundary.
Since $W\times T$ is a deformation retract of $W^+\times T^+$, we can extend
$H$ to an embedding $H:W^+\times T^+/
(x\times t\sim x\times\sigma\mid x\in\ofk,t\in T)\cong W^+\times T^+
\hookrightarrow M$ such that
for this extended $H$ we still have $H|_{W^+\times\sigma_i}$ --- a smooth
embedding into $f_i(V_i)$, disjoint from self-intersections of $f_i(V_i)$,
and $H(W^+\times T^+-\ofk\times T-W^+\times\{\sigma_1,\dots,\sigma_k\})
\cap(\bigcup_{i=1\dots k}f_i(V_i))=\emptyset$.

Denote by $\pi_W:W^+\times T^+\to W^+$ the natural projection. By definition,
$\nu_{W\times T}\simeq
\pi_W^*(J^*(\pi_1^*\nu_{f_1}\oplus\dots\oplus\pi_k^*\nu_{f_k})
\oplus J^*\Phi(\cdot,\sigma)^*\nu_M)$. On the other hand,
$\nu_{W\times T}=\nu_H\oplus H^*\nu_M$. It is easy to see that
$H$ is homotopic to $\Phi(\cdot,\sigma)\circ J\circ\pi_W$. Therefore
$\nu_H\simeq\pi_W^* J^*(\pi_1^*\nu_{f_1}\oplus\dots\oplus\pi_k^*\nu_{f_k})$.
Since $\dim(W\times T)=p+k<m-n_i=\dim(\pi_W^*J^*\pi_i^*\nu_{f_i})$, then
$\pi_W^*J^*\pi_i^*\nu_{f_i}=\eta_i\oplus\epsilon_i$ for some bundle $\eta_i$
and a trivial line bundle $\epsilon_i$. So,
$\nu_H\simeq\eta_1\oplus\dots\oplus\eta_k$. Since
$\dim(\nu_H)=m-p-k=\dim(\eta_1\oplus\dots\oplus\eta_k)>p+k=\dim(W\times T)$,
then $\nu_H\cong\eta_1\oplus\dots\oplus\eta_k$. Fix such an isomorphism.
Denote by $\delta_i$ the
barycenter of the face, opposite to $\sigma_i$ in the simplex $T$. Denote by
$\varepsilon_i$ the trivial 1-dimensional bundle over $W^+\times T^+$,
parallel to the line $\sigma_i\sigma\delta_i\subset T$. Note that
$\varepsilon_1,\dots,\varepsilon_k$ are linearly dependent
($\vec{\varepsilon}_1+\dots+\vec{\varepsilon}_k=0$), but any $(k-1)$ of
them are linearly independent. Denote by $\Xi$ the bundle, spanned by
$\varepsilon_1,\dots,\varepsilon_k$, then $\Xi$ is the trivial bundle,
tangent to $T^+$ in $W^+\times T^+$. Recall that
$H(W^+\times\sigma_i)\subset f_i(V_i)$. Therefore
$\nu_{f_i}(V_i,M)|_{H(W^+\times\sigma_i)}\subset
\nu_{H|_{W^+\times\sigma_i}}=\nu_H(W^+\times T^+,M)|_{W^+\times\sigma_i}\oplus
\nu(W^+\times\sigma_i,W^+\times T^+)$. In
$\nu_H|_{W^+\times\sigma_i}\oplus\nu(W^+\times\sigma_i,W^+\times T^+)$
we already have a summand, isomorphic to the pull-back of $\nu_{f_i}$,
namely $\eta_i\oplus\varepsilon_i$. Therefore the complement to
$\nu_{f_i}(V,M)$ in $\nu_{H|_{W^+\times\sigma_i}}$ is stably isomorphic
to the complement of $\eta_i\oplus\varepsilon_i$. Since
$\dim(\nu_{H|_{W^+\times\sigma_i}})-\dim(\nu_{f_i}(V_i,M))=n_i-p-1>p+1=
\dim(W^+\times\sigma_i)$, then these complements are isomorphic.
So, in $\nu_H\oplus\Xi$ we have $k$ bundles
$\chi_i=(\eta_1\oplus\dots\oplus\widehat{\eta_i}\oplus\dots\oplus\eta_k)
\oplus\varepsilon_i^\perp(\Xi)$, and $\chi_i|_{W^+\times\sigma_i}$
is isomorphic to the complement to $\nu_{f_i}(V_i,M)$ in
$\nu_{H|_{W^+\times\sigma_i}}=(\nu_H\oplus\Xi)|_{W^+\times\sigma_i}$. Here
$\varepsilon_i^\perp(\Xi)$ is the orthogonal complement to $\varepsilon_i$
in~$\Xi$, then $\varepsilon_i^\perp(\Xi)|_{W^+\times\sigma_i\sigma\delta_i}=
\nu(W^+\times\sigma_i\sigma\delta_i,W^+\times T^+)$. Note that common
intersection of $\chi_i|_{W^+\times\sigma}$ is empty. Since the homotopy group
of the Stiefel manifold $\pi_d(V_{\dim(\nu_H\oplus\Xi),\dim(\chi_i)})=0$ if
$d<\dim(\nu_H\oplus\Xi)-\dim(\chi_i)$ \cite{Whd}, and
$\dim(W\times T)=p+k<m-n_i=\dim(\nu_H\oplus\Xi)-\dim(\chi_i)$, then there
is no obstruction for homotopying $\chi_i$ in a small neighborhood of
$W^+\times\sigma_i$ in $W^+\times T^+$ to a bundle $\chi_i'$ such that
$\chi_i'|_{W^+\times\sigma_i}=
\nu_{f_i}(V_i,M)^\perp(\nu_{H|_{W^+\times\sigma_i}})$.

Let $\varphi:W^+\to\R$ be a smooth function such that $\varphi^{-1}(0)=\ofk$,
$\varphi^{-1}(1)=N$ and $\varphi^{-1}[0,1]=W$. Since $H$ is an embedding,
then on each sheet $H(W^+\times\delta_i\sigma\sigma_i^+)$ we can introduce
a coordinate system $W^+\times\R$, such that $H(W^+\times\sigma)$ have
coordinates $(w,0)$, and $H(W^+\times\sigma_i)$ have coordinates
$(w,\varphi(w))$. Let $U_i$ be a small enough neighborhood of the zero-section
of $\chi_i|_{W^+\times\delta_i\sigma\sigma_i^+}$ such that
$\exp\circ dH:U_i\to M$ is a diffeomorphism on its image. On the manifold
$\exp\circ dH(U_i)$ we have coordinate system $(w,t,\vec{v})$, where
$w\in W,t\in\R,\vec{v}\in U_i$. Since $\chi_i|_{W^+\times\sigma_i}=
\nu_{f_i}(V_i,M)^\perp(\nu_{H|_{W^+\times\sigma_i}})$, then
$\exp\circ dH(U_i|_{W^+\times\sigma_i})$ is a tubular neighborhood of
$H(W^+\times\sigma_i)$ in $f_i(V_i)$. Let $\psi:W^+\to[0,1]$ be a smooth
function such that $\psi\geq0$, $\psi|_{W}\equiv1$ and $\psi|_{W^+-U}\equiv0$,
where open set $U$ is such that $W\subset U\subset W^+$ and
$\overline{U}\subset W^+$. Let $\kappa:\R\to\R$ be a smooth bell-shaped
function such that $\kappa(0)=1$ and $\kappa(\|\vec{v}\|)=0$ for $\vec{v}$
outside of $U_i$. Then the regular homotopy from $f_i$ to $f_i'$ is given by
$(w,\varphi(w),\vec{v})\mapsto
(w,\varphi(w)-t\psi(w)\kappa(\|\vec{v}\|),\vec{v})$, $t\in[0,1]$ for points
in $\exp\circ dH(U_i|_{W^+\times\sigma_i})\subset f_i(V_i)$, and is constant
outside. Since common intersection of $\chi_i|_{W^+\times\sigma}$ is empty,
then common intersection of $\exp\circ dH(U_i)$ is empty, and at the last
moment we will get exactly $\op(f_1'\times\dots\times f_k')=N$.
\end{proof}
\begin{proof}[Proof of Theorem~\ref{egroup}]
If $\wF=\emptyset$, put $b(F)=0\in\Omega_p(E_k;\xi(f))$. Now suppose
$\wF\neq\emptyset$. Denote by $D:M\to\DM$ the diffeomorphism
$x\mapsto(x,\dots,x)$. Obviously, the restriction on $\wF$ of the projection
$\Vk\to V$ on the $i$-th factor is a composite $h_i\pi_i'$ for a uniquely
determined $h_i:\wF_{(i)}\to V$, and $D^{-1}F|_\wF$ is a composite $h\pi'$
for a uniquely determined $h:\op(F)\to M$. Since $\wF\subset(\Vk-\DV)$, then
$h_i\pi_i'(x)\neq h_j\pi_j'(x)$ for $x\in\wF$ and $i\neq j$. Let
$\bar F:\Vk\times[0,1]\to\Mk$ be a $\sk$-equivariant homotopy between
$\bar F_0=\fk$ and $\bar F_1=F$. Denote by $\pi^M_i:\Mk\to M$ the projection
on the $i$-th factor. Then $fh_i\pi_i'$ and $D^{-1}F=\pi^M_iF$ are connected
by the homotopy $x\mapsto\pi^M_i\bar F_t(h_1\pi_1'(x),\dots,h_k\pi_k'(x)),
t\in[0,1]$. This gives a $\sk$-invariant homotopy $H:\wF\times T^1\to M$,
where $T^1$ is the $1$-dimensional subcomplex of $T$, consisting of vertices
$\sigma,\sigma_1,\dots,\sigma_k$ and straight intervals
$\sigma\sigma_1,\dots,\sigma\sigma_k$. Since $T^1$ is an $\sk$-invariant
deformation retract of $T$ under $\sk$-equivariant deformation, we can extend
$H$ to a $\sk$-invariant homotopy $H:\wF\times T\to M$. Note that if $F=\fk$,
we can take $H$ to be the constant homotopy. By universal property of
$\widehat{E_k}$, this gives a canonical map $j:\op(F)\to E_k$.

Since the manifold $\wF\hookrightarrow\Vk$, then
$\nu_\wF=\nu_\Vk|_\wF\oplus\nu(\wF,\Vk)$. By construction,
$\nu(\wF,\Vk)\cong(F|_\wF)^*\nu(\DM,\Mk)$ and
$\nu_\Vk|_\wF\cong\nu_\fk|_\wF\oplus(\fk|_\wF)^*\nu_\Mk$. Since $\fk$ is
$\sk$-equivariantly homotopic to $F$ and $F(\wF)\subset\DM$, then
$\nu_\wF\cong\nu_\fk|_\wF\oplus(F|_\wF)^*(\nu_\Mk|_\DM\oplus\nu(\DM,\Mk))$.
Note that $F|_\wF=D\circ\Phi(\cdot,\sigma)\circ\wj$ and $\fk|_\wF=
f^\k\circ(\pi_1\times\dots\times\pi_k)\circ\wj$. Since
$\nu_\Mk|_\DM\oplus\nu(\DM,\Mk)\simeq\nu_\DM$ and $\Phi(\cdot,\sigma)$ is
$\sk$-invariant, then $\nu_\wF\simeq\wj^*(\pi_1\times\dots\times\pi_k)^*
\nu_\fk\oplus\wj^*\Phi(\cdot,\sigma)^*\nu_M$. Therefore $\nu_{\op(F)}=
\nu_\wF/\sk\simeq j^*(\pi_1^*\nu_f\oplus\dots\oplus\pi_k^*\nu_f)/\sk\oplus
j^*\phi^*\nu_M=j^*\xi$. Thus the map $j:\op(F)\to E_k$ defines
an element $b(F)\in\Omega_p(E_k;\xi(f))$.

Now, suppose we have a $k$-disjoint homotopy $\bar F:\Vk\times[1,2]\to\Mk$,
where $F_1=\bar F_1,F_2=\bar F_2$ are $k$-transversal. Since $\bar F$ is
$k$-disjoint, then $\bar F^{-1}(\DM)-\DV\times[1,2]$ is disjoint from
$\DV\times[1,2]$. WLOG we may assume $\widehat{\op(F_1)}\neq\emptyset$. If
$\widehat{\op(F_2)}\neq\emptyset$, choose an open $\sk$-invariant neighborhood
$U$ of $\DV\times[1,2]$ in $\Vk\times[1,2]$ such that $\overline{U}$ is
disjoint from $\bar F^{-1}(\DM)-\DV\times[1,2]$.
If $\widehat{\op(F_2)}=\emptyset$, we may assume that
$(\bar F_t)^{-1}(\DM)-\DV=\emptyset$ for $t$ close to~$2$. Then take $U$ to
be an open $\sk$-invariant neighborhood of $\DV\times[1,2]\cup\Vk\times\{2\}$
in $\Vk\times[1,2]$ such that $\overline{U}$ is disjoint from
$F^{-1}(\DM)-\DV\times[1,2]$.

Note that the $\sk$-action on $\Vk\times[1,2]-U$ is free. Therefore we can
approximate $\bar F|_{\Vk\times[1,2]-U}$ by a smooth $\sk$-equivariant map
$\tilde F:\Vk\times[1,2]-U\to\Mk$, which is transversal to $\DM$, keeping
it fixed on the end $\Vk\times\{1\}-U$, where it is already transversal
(and on $\Vk\times\{2\}-U$, if $\widehat{\op(F_2)}\neq\emptyset$). Since
$\bar F(\partial U)\cap\DM=\emptyset$, then $\tilde F(\partial U)\cap\DM=
\emptyset$, and $\wW=\tilde F^{-1}(\DM)$ is a proper compact submanifold of
$\Vk\times[1,2]-U$ with a free $\sk$-action, and
$\partial\wW=\widehat{\op(F_1)}\cup\widehat{\op(F_2)}$. Denote $W=\wW/\sk$.

Obviously, the restriction on $\wW$ of the projection $\Vk\times[1,2]\to V$
on the $i$-th factor is a composite $h_i\pi_i'$ for a uniquely determined
$h_i:\wW_{(i)}\to V$, and $D^{-1}\tilde F|_\wW$ is a composite $h\pi'$ for
a uniquely determined $h:W\to M$. Since $\wW\subset(\Vk-\DV)\times[1,2]$,
then $h_i\pi_i'(x)\neq h_j\pi_j'(x)$ for $x\in\wW$ and $i\neq j$. Let
$\tilde F:\Vk\times[0,1]\to\Mk$ be a $\sk$-equivariant homotopy between
$\tilde F_0=\fk$ and $\tilde F_1=\bar F_1=F_1$. Denote by $t:\Vk\times[0,2]
\to[0,2]$ the projection on the last factor. Then $fh_i\pi_i'$ and
$D^{-1}(\tilde F|_\wW)=\pi^M_i(\tilde F|_\wW)$ are connected by the homotopy
$x\mapsto\pi^M_i\tilde F(h_1\pi_1'(x),\dots,h_k\pi_k'(x),t't(x)),t'\in[0,1]$.
This gives a $\sk$-invariant homotopy $H:\wW\times T^1\to M$. Since $T^1$
is an $\sk$-invariant deformation retract of $T$ under $\sk$-equivariant
deformation, we can extend $H$ to a $\sk$-invariant homotopy
$H:\wW\times T\to M$. By universal property of $\widehat{E_k}$, this gives
a canonical map $J:W\to E_k$.

Since the manifold $\wW\hookrightarrow\Vk\times[1,2]$, then
$\nu_\wW=\nu_{\Vk\times[1,2]}|_\wW\oplus\nu(\wW,\Vk\times[1,2])$. By
construction, $\nu(\wW,\Vk\times[1,2])\cong(\tilde F|_\wW)^*\nu(\DM,\Mk)$
and $\nu_{\Vk\times[1,2]}|_\wW\cong(\pi_1^V\times\dots\times\pi_k^V|_\wW)^*
(\nu_\fk\oplus(\fk)^*\nu_\Mk)$, where $\pi_i^V:\Vk\times[0,2]\to V$ is the
projection on the $i$-th factor. Since
$\fk\circ(\pi_1^V\times\dots\times\pi_k^V)|_\wW$ is $\sk$-equivariantly
homotopic to $\tilde F|_\wW$ and $\tilde F(\wW)\subset\DM$, then
$\nu_\wW\cong(\pi_1^V\times\dots\times\pi_k^V|_\wW)^*\nu_\fk\oplus
(\tilde F|_\wW)^*(\nu_\Mk|_\DM\oplus\nu(\DM,\Mk))$. Note that
$\tilde F|_\wW=D\circ\Phi(\cdot,\sigma)\circ\wJ$
and $\fk\circ(\pi_1^V\times\dots\times\pi_k^V)|_\wW=
f^\k\circ(\pi_1\times\dots\times\pi_k)\circ\wJ$. Since
$\nu_\Mk|_\DM\oplus\nu(\DM,\Mk)\simeq\nu_\DM$ and $\Phi(\cdot,\sigma)$ is
$\sk$-invariant,
then $\nu_\wW\simeq\wJ^*(\pi_1\times\dots\times\pi_k)^*\nu_\fk\oplus
\wJ^*\Phi(\cdot,\sigma)^*\nu_M$. Therefore $\nu_W=\nu_\wW/\sk\simeq
J^*(\pi_1^*\nu_f\oplus\dots\oplus\pi_k^*\nu_f)/\sk\oplus J^*\phi^*\nu_M=
J^*\xi$. Thus the map $J:W\to E_k$ gives a bordism between $b(F_1)$ and
$b(F_2)$ in $\Omega_p(E_k;\xi(f))$.
\end{proof}
\begin{proof}[Proof of Theorem~\ref{center}]
Making a small regular perturbation of $f$, we may assume that it is in
general position. By Theorem~\ref{egroup}, this will not change the class
$b(\fk)\in\Omega_p(E_k;\xi)$. From $(k+1)(n+1)\leq km$ it follows that $2p<n$,
or $2\dim(\op(k,f))<\dim(V)$. Since $\widehat{\op(\fk)}\subset\Vk-\DV$, by
general position we may assume that the manifold $\widehat{\op(\fk)}_{(i)}$
is embedded by $h_i$ into $V$ for all $i=1\dots k$, the manifold $\op(\fk)$
is embedded by $h$ into $M$, and
$f:\left(\bigcup_ih_i(\widehat{\op(\fk)}_{(i)})\right)\to h(\op(\fk))$ is a
$k$-fold covering (possibly, some of $h_i(\widehat{\op(\fk)}_{(i)})$ coincide).

Let $J:W\to E_k$ be the bordism between $\op(\fk)$ and $N$ in
$\Omega_p(E_k;\xi)$. Choose a lifting $\wJ:\wW\to\widehat{E_k}$. Since the
$\sk$-action on $\wW$ is free, the diagonal $\sk$-action on $\wW\times T$ is
also free. Therefore $(\wW\times T)/\sk$ is a smooth manifold. Note that
$(\wW\times\sigma)/\sk=W$. Since $\Phi:\widehat{E_k}\times T\to M$ is
$\sk$-invariant, the map $B:(\wW\times T)/\sk\to M$ is well-defined by the
formula $[w,t]\mapsto\Phi(\wJ(w),t)$. Note that
$B((\wW\times\sigma_i)/\sk)\subset f(V)$. By construction of $b(\fk)$, for
$x\in\widehat{\op(\fk)}$ and $t,t'\in T$ we have $B([x,t])=B([x,t'])$.

Since $2(p+1)<n$ and $(p+1)+(2n-m)<n$, we have $2\dim(W)<\dim(V)$, and
$\dim(W)+\dim($ self-intersections of $V$ in $M)<\dim(V)$. Applying general
position argument, we may $C^0$ perturb the map $B:(\wW\times T)/\sk\to M$
(leaving it fixed on $(\widehat{\op(\fk)}\times T)/\sk$) so that
$B|_{(\wW\times\sigma_i)/\sk}$ will be a smooth embedding into $f(V)$,
disjoint from self-intersections of $f:V\looparrowright M$, distinct from
$\op(\fk)$.

Since $2(p+k)<m$ and $(p+k)+n<m$, we have $2\dim(W\times T)<\dim(M)$,
and $\dim(W\times T)+\dim(V)<\dim(M)$. This means that by general position
we can $C^0$-perturb the map $B$ (leaving it fixed on
$(\widehat{\op(\fk)}\times T)/\sk$ and $(\wW\times\sigma_i)/\sk$)
so that $B:\left(\wW\times T-\widehat{\op(\fk)}\times T\right)/\sk
\hookrightarrow M$ will be a smooth embedding with
$B\left(\left\{\wW\times T-\widehat{\op(\fk)}\times T-
\wW\times\{\sigma_1,\dots,\sigma_k\}\right\}/\sk\right)\cap f(V)=\emptyset$.

Denote by $\wW^+$ and $T^+$ open manifolds without boundary, which
are obtained from $\wW$ and $T$ by attaching a small collar neighborhood
of the boundary. Then $\left(\wW^+\times T^+/\left\{
x\times t\sim x\times\sigma\mid x\in\widehat{\op(\fk)},t\in T\right\}
\right)/\sk\cong(\wW^+\times T^+)/\sk$ is a smooth manifold. Since
$\wW\times T$ is a $\sk$-equivariant deformation retract of $\wW^+\times T^+$,
we can extend $B$ to an embedding $B:\left(\wW^+\times T^+/\left\{
x\times t\sim x\times\sigma\mid x\in\widehat{\op(\fk)},t\in T\right\}
\right)/\sk\hookrightarrow M$ such that for this extended $B$ we still
have $B|_{(\wW^+\times\sigma_i)/\sk}$ --- a smooth embedding into $f(V)$,
disjoint from self-intersections of $f(V)$, and
$B\left(\left\{\wW^+\times T^+-\widehat{\op(\fk)}\times T-
\wW^+\times\{\sigma_1,\dots,\sigma_k\}\right\}/\sk\right)\cap
f(V)=\emptyset$.

Denote by $H:\wW^+\times T^+\looparrowright M$ the immersion, defined
as a composition of the natural covering
$\wW^+\times T^+\to(\wW^+\times T^+)/\sk$ with $B$. Strictly speaking,
this is not quite immersion, it has "singularities" at $\op(\fk)\times T$.
Obviously, $H$ is $\sk$-invariant. The rest of the proof follows the proof
of Theorem~\ref{easy}, which is expressly written in a $\sk$-invariant
language. Indeed, substitute $W\mapsto\wW$, $\nu_{f_i}\mapsto\nu_f$,
$V_i\mapsto V$, $J\mapsto\wJ$. Locally, in a small neighborhood of $w\in W$,
the construction from the proof of Theorem~\ref{easy} translates without
changes. Then the fact that $H$ is $\sk$-invariant guarantee us that the
desired regular homotopy will be well-defined globally.
\end{proof}

\end{document}